\documentclass[11pt , reqno]{article}
\usepackage{amssymb,amsmath,amsfonts,amsthm}
\usepackage{enumerate}
\usepackage{graphicx}
\usepackage[numbers]{natbib}
\usepackage{float}
\makeatletter
\DeclareRobustCommand*{\bfseries}{%
  \not@math@alphabet\bfseries\mathbf
  \fontseries\bfdefault\selectfont
  \boldmath
} \makeatother

\newtheorem{thm}{Theorem}[section]

\begin{document}

\title{\bf A  generalization of a theorem of Nash-Williams}

\date{\vspace{-1.5cm} Compiled on \today \vspace{1cm}}       
\maketitle

\vspace{-50pt}

\noindent $\mathrm{D.\, Bauer}^{\mathrm{a}},\, \mathrm{L.\, Lesniak}^{\mathrm{b}}, \, \mathrm{E.\, Schmeichel}^{\mathrm{c}}$
 \\
{\footnotesize  $^{\mathrm{a}}$ \textit{Department of Mathematical
Sciences,
Stevens Institute of Technology, Hoboken, NJ 07030, USA, dbauer@stevens.edu} \\
\footnotesize $^{\mathrm{b}}$ \textit{Department of Mathematics, Western Michigan University, Kalamazoo, MI 49008, USA,  \\linda.lesniak@wmich.edu} \\
\footnotesize $^{\mathrm{c}}$ \textit{Department of Mathematics,
San Jos\'{e} State University, San Jos\'{e}, CA 95192, USA, \\edward.schmeichel@sjsu.edu} \\

\begin{abstract}

In \cite{Chvatal72}, Chv\'{a}tal  gave a well-known sufficient condition for a graphical sequence to be forcibly hamiltonian, and showed that in some sense his condition is best possible.  Nash-Williams  \cite{NW71}  gave examples of forcibly hamiltonian $n-$sequences that do not satisfy Chv\'{a}tal's condition, for every  $n\geq 5$.  In this note we generalize the Nash-Williams examples, and use this generalization to generate $\Omega ( \frac{2^n}{\sqrt n} )$ forcibly hamiltonian $n-$sequences that do not satisfy Chv\'{a}tal's condition.
\end{abstract}

\section{Introduction}

Our terminology and notation are standard except as indicated.  A good reference for undefined terms is \cite{CLZ16}.

A graphical sequence $\pi = (d_1 \leq ... \leq d_n) $ is said to be \emph{forcibly hamiltonian} if every graph with degree sequence $\pi$ is hamiltonian.  Beginning in 1952, a series of conditions for $\pi$ to be forcibly hamiltonian was given by Dirac \cite{Dirac}, P\'{o}sa \cite{Posa}, and Bondy \cite{Bondy}, culminating in the following condition of Chv\'{a}tal \cite{Chvatal72} in 1972.

\begin{thm}\label{thm:COT}  Let $\pi = (d_1 \leq \cdots \leq d_n)$ be a graphical sequence with $n \geq 3$.  If 
	
		\begin{equation}
	d_k \geq k+1  \text{ or }   d_{n-k} \geq n-k, \quad \text{ for } 1 \leq k \leq \frac{n-1}{2},
	\end{equation}
	
	  then $\pi$ is forcibly hamiltonian.
\end{thm}

Chv\'{a}tal  also showed that his condition for a graphical sequence to be forcibly hamiltonian is best possible in the following sense.  Suppose a graphical  sequence $\pi = (d_1 \leq \cdots \leq d_n)$ fails to satisfy Chv\'{a}tal's condition.  Then $d_k \leq k$ and $d_{n-k} \leq n - k -1$, for some $1 \leq \ k \leq (n-1)/2$.  Given positive integers $k,n$ for which $1 \leq k \leq (n-1)/2$, let $C_{n,k}$ denote the graph $K_k \vee (\bar{K_k} \cup K_{n-2k})$.  Clearly $C_{n,k}$ is not hamiltonian.  Moreover, this graph has degree sequence $\pi^* = (d_1^*  \leq \cdots \leq d_n^* )$, where $d_1^* = \cdots = d_k^*  = k,$    $d_{k+1}^* = \cdots = d_{n-k}^* = n-k-1,$ and $d_{n-k+1}^* = \cdots = d_n^* = n-1.$  Finally, it is immediate that $d_j^* \geq d_j$ for $ 1 \leq j \leq n$.

Of course, not every forcibly hamiltonian sequence $\pi = (d_1 \leq \cdots \leq d_n) $ satisfies Chv\'{a}tal's condition (1).  Perhaps the best-known examples are the $(n-1)/2-$regular $n-$sequences of Nash-Williams , where $n\equiv 1$ (mod 4) and $n \ge 5.$  More generally, Nash-Williams (see \cite{NW71})  gave the following families of forcibly hamiltonian sequences that do not satisfy (1).

\begin{thm}\label{thm:NW}  For each  integer $n \geq 5$ and even positive integer $k\leq  \frac{n-1}{2}$, the graphical sequence $\pi = (d_1  \leq \cdots  \leq d_n)$ with $d_j = k$ for $j = 1,  \cdots, k,$ and $d_j = n-k-1$ otherwise, is forcibly hamiltonian. 

\end{thm}

In this note we give a generalization of Theorem 1.2.   In particular, we prove that the final $k-1$ terms of $\pi$ in Theorem 1.2 can be altered essentially arbitrarily without affecting the conclusion.  We then show how to use this generaliztion to produce exponentially many graphical sequences that are forcibly hamiltonian, but do not satisfy Chv\'{a}tal's condition (1).

\section{Generalization of Nash-Williams' Theorem 1.2}

The proof of our generalization of Theorem 1.2 depends on two preliminary results.  

First, the following  result of Dirac \cite{Dirac} gives a  lower bound on the circumference $cir(H)$ of a $2$-connected graph $H$ in terms of its order and minimum degree.

\begin{thm} If $H$ is a 2-connected graph of order $n$ with minimum degree $\delta(H)$, then  the circumference $cir(H)$ satisfies

$$cir(H) \geq  min \{ n, 2\delta(H) \}.$$

\end{thm}

Second, the \emph{closure} $cl(G)$ of a graph $G$ of order $n$ is the graph obtained from $G$ by recursively joining pairs of vertices whose degree sum is at least $n$ (in the resulting graph at each stage) until no such pair exists.  Bondy and Chv\'{a}tal \cite{BC}  observed  the following important result relating the hamiltonicity of a graph and its closure.

\begin{thm}\label{thm:BC} A graph is hamiltonian if and only if its closure is hamiltonian.

\end{thm}

Theorem 2.3 gives our generalization of  Theorem 1.2.

\begin{thm}\label{thm:BLS}  Let $n \geq 5$ and $k \geq 2$ be integers with $k < n/2$, and let $\pi = (d_1 \leq \cdots  \leq d_n) $ be a graphical sequence. If $d_1 = \cdots = d_k  = k$ and $ d_{k+1} = \cdots = d_{n-k+1} = n-k-1$, then $\pi$ is forcibly hamiltonian unless  

\begin{equation}
	d_j = n-k-1 \text { for } n-k+2  \leq  j  \leq  n-1  \text { and } d_n = n-1. 
	\end{equation}

\end{thm}

\begin{proof}  Assume first that $n$ is odd and $k = \frac{n-1}{2}.$  Then $n-k+1 = \frac{n+3}{2}$ and $n-k-1 = \frac{n-1}{2}$, and so $d_1 = \cdots = d_{(n+3)/2} = \frac{n-1}{2}.$

 Let $G$ be any graph with degree sequence $\pi.$   We claim that $G$ is 2-connected or (2) occurs.  Suppose first $G$ is disconnected, with $G = H_1 \cup H_2.$  Then $d_1 = \frac{n-1}{2}$ implies that $|H_i|\geq \frac{n+1}{2}$ for each $i$, producing a contradiction since $G$ has order $n$.  Similarly, if $G$ has a cut-vertex $v$ and $G - v =  H_1 \cup H_2$, then $d_1 = \frac{n-1}{2}$ implies that $|H_1| = |H_2| = \frac{n-1}{2}$, and that every vertex in $H_1 \cup H_2$  is adjacent to $v$.  Thus, $d_j = \frac{n-1}{2}$ for $j = 1 \cdots n-1$ and $d_n = n-1,$ and (2) occurs.

Thus we may assume that $G$ is 2-connected and may apply Theorem 2.1 to conclude that  $cir(G) \geq  min \{ n, 2d_1 \} \geq n-1.$  Let $C$ be a longest cycle of $G$.  We claim that $C$ is a hamiltonian cycle.  Assume, to the contrary, that $C$ has length $n-1,$ and that  $v$ is the vertex of $G$ not on $C$.  Since $C$ is a longest cycle, $v$ cannot be adjacent to consecutive vertices of $C$.  Then $d_1 = \frac{n-1}{2}$ implies that $d(v) = \frac{n-1}{2}, $ and $v$ is adjacent to alternate vertices of $C$.

If two vertices of $I := V (G-v) - N_G(v)$ are adjacent, then $G$ is hamiltonian.  Thus $I$ is an independent set.  Since every vertex in I has degree at least $d_1=\frac{n-1}{2}$, it follows that every vertex in $I$ is adjacent to every vertex in $N(v)$, so that every vertex in $N(v)$ has degree at least $|I \cup \{v\}| = \frac{n+1}{2}.$  Since $|N(v)| = \frac{n-1}{2},$ this contradicts the assumption that $d_{(n+3)/2} = \frac{n-1}{2}.$  Thus the proof is complete when  $n$ is odd and $k = \frac{n-1}{2}.$

Therefore,  we may assume that $k \leq \frac{n}{2} - 1,$ so that $n-k-1 \geq \frac{n}{2}.$   Let $G$ be any graph with degree sequence $\pi.$   Let $A$ denote the vertices of $G$ of degree $k$, define $B$ to be $N(A) - A$ and let $C = V(G) - A - B.$  Since $n-k-1 \geq \frac{n}{2},$ it follows that $B \cup C$ induces a complete graph in the closure $cl(G).$  Since every vertex $v$ in $ B$ is also adjacent to a vertex in $A,$ it follows that each vertex $v$ in $ B$ has degree at least $n-k$ in $cl(G),$  and thus each such  $v $ in $B$ is adjacent to every vertex of $A$ in $cl(G).$ 

 Finally, we claim that there are no edges between $A$ and $C$ in $cl(G).$  Suppose, to the contrary, that $xy$ is the first edge between vertices of $A$ and $C$ added in forming $cl(G),$ where $x \in A$ and $y \in C.$  When $xy$ is added to $cl(G),$ we have $d_{cl(G)}(y) = n-k-1$ and so $d_{cl(G)}(x) \geq k+1.$  But then $x$ will be adjacent to \emph{every} vertex of $C$ in $cl(G),$   implying that every vertex of $C$ has degree at least $n-k$ in $cl(G).$  Since every vertex of $A$ has degree at least $k$ in $cl(G),$ it follows that every vertex of $A$  will be adjacent to every vertex of $C$ in $cl(G).$  Consequently, $cl(G)$ contains $\bar{K_k} \vee K_{n-k},$ where $n-k > k,$ and so $cl(G)$, and consequently $G,$  is hamiltonian. So we may assume that no vertices of $A$ and $C$ are adjacent in $cl(G)$ and, in particular, each vertex of $A$ has degree $k$ in $cl(G),$ as well as in $G.$

Setting $|B| = j,$ it follows that  $K_j \vee (\bar{K_k} \cup K_{n-j-k})$ is a spanning subgraph of $cl(G).$  Consequently, $cl(G) $is hamiltonian provided $j > k.$   Thus we may assume that $j \leq k.$  

Since each vertex in $A$ has degree  $k$ in $cl(G),$ the subgraph $H$ induced by $A$ in $cl(G)$ is $(k-j)-$ regular. Consequently, $H$ has a path  containing at least $k-j+1$ vertices.  If, in fact, any longest path $P$ in $H$ contains at least $k-j+2$ vertices, then the vertex set of $H$ can be partitioned into the vertex set of the path $P$ plus the remaining (at most $j-2$)  vertices.  Therefore, $cl(G)$ is hamiltonian.  Assume then that every longest  path in $H$ contains exactly $k-j+1$ vertices and let $P : v_0, v_1, \cdots, v_{k-j}$ be such a path.  Then $v_0$ is adjacent in $H$ only  to the $k-j$ other vertices of $P$.  But since $v_0,  v_, \cdots, v_{k-j}, v_0$ is  a cycle, this is true of every vertex of $P$.  Thus the $k-j+1$ vertices of $P$ induce a complete graph in $H.$   Iterating this process,we can partition $H$ into  $\frac{k}{k-j+1}$ disjoint complete graphs of order $k-j+1.$

Since each of the $\frac{k}{k-j+1}$ components of $H$ is complete (and hence has a spanning path), if $\frac{k}{k-j+1} \leq j-1,$ then, as before $cl(G)$ is hamiltonian.  On the other hand, since $1 \leq j \leq k,$  it follows that $\frac{k}{k-j+1} \geq j$ only for $j = 1$ or $j=k$, where equality occurs in each case. We conclude that for $cl(G)$ to be nonhamiltonian, either $j=1$ and $cl(G)= K_1 \vee (K_k \cup K_{n-k-1} ),$ or $j=k$ and $cl(G)= K_k \vee (\bar{K_k} \cup K_{n-2k}) = C_{n,k}$. 

We finally consider  the graph $G$ itself giving rise to such a closure $cl(G).$

For $j=1,$ it is clear that for the vertices of $A$ to have degree $k$ in $G,$ the subgraph induced by $A \cup B$ in $G$ must be complete.   Similarly, for the vertices in $C$ to have degree
 at least $n-k-1$ in $G,$ the subgraph induced by $B \cup C$ in $G$ must be complete.  But then the single vertex of $B$ has degree $n-1$  in $G$ and the degrees of $G$ satisfy (2).

For $j=k,$  noting that each vertex in $A$ (repsectively, $C$) has degree $k$  (respectively, at least $n-k-1$) in $G,$ it is clear that every vertex  in both $A$ and $C$ must be adacent to every vertex of $B$ in $G.$  So every vertex of $B$ has degree at least $|A \cup C| = k + (n-2k) = n-k > n-k-1$ in $G.$ On the other hand, $\pi$ contains at least $n-2k+1$ terms equal to $n-k-1,$ and since $|C|=n-2k,$ at least one vertex of $B$ must have degree $n-k-1$ in $G.$  This contradiction completes the proof.

\end{proof}

We note that any graphical sequence that is forcibly hamiltonian by Theorem 2.3 fails Chv\'{a}tal's condition (1) for the value of $k$ described in the theorem.  We call such a sequence a \emph{Nash-Williams (n,k)-sequence.}  In the next section we construct new families of Nash-Williams sequences.

\section{Constructing Nash-Williams Sequences}

In this section we apply Therem 2.3 to construct families of Nash-Williams sequences.  In order to do so, we will use a version of the well-known theorem of Erd\"{o}s and Gallai \cite{EG} to show that an $n-$sequence is graphical.  Note that in order to apply this theorem, the sequence under consideration must be in nonincreasing order.   The form of this theorem (see \cite{SH91}) we will use depends on the value $D= $ max  $\{i : d_i \geq i\}.$

\begin{thm} A sequence  $(d_1 \geq \cdots \geq d_n)$, with $n \geq 2$, of nonnegative integers is graphical if and only if the sequence has even sum and

\begin{equation}
	\sum_{i=1}^{j} d_i \leq j(j-1)   + \sum_{i=j+1}^{n} \text{ min } \{j, d_i\},  \text{ for }  d_n <  j  \leq D.
	\end{equation}

\end{thm}

Given integers $n \geq 5 $  and $k \geq 2$, with $k < n/2,$  we define two \emph{foundational $(n,k)-$sequences}  $(s_1  \leq \cdots \leq s_n)$ as follows.  In both sequences,  $s_j = k$  for    $\j = 1, \cdots, k,$ and $s_j = n-j-1$  for $j = k+1, \cdots, n-1.$  The final term is $n-k-1$ in the first sequence and $n-k$ in the second sequence.  We note that the first foundational $(n,k)-$sequence has even sum if and only if $k$ is even while  the second foundational $(n,k)-$sequence has odd sum if and only if $k$ is even.

Let $\pi' = (d_1' \leq \cdots \leq d_{k-1}')$ be any  sequence of integers with  $ 0 \leq d_j' \leq k-1$ and $\pi' \neq (0, \cdots  , 0 , k-1).$  Given $\pi'$, select the foundational $(n,k)-$sequence $s$ with the same parity of sum as $\pi',$ and then form the $n-$sequence $\pi$ by adding in order the terms of $\pi'$ to the final $k-1$ terms in $s.$  We claim that $\pi$ is forcibly hamiltonian.

If $\pi$ is graphical, then since $\pi$ satisfies the conditions in Theorem 2.3 (in particular,  $\pi' \neq (0, \cdots , 0, k-1)$ implies that $\pi$  is not the sequence described in (2)) it follows that $\pi$ is forcibly hamiltonian.  Since, by our construction, $\pi$ has even sum,  to show that $\pi$ is graphical we need only  show that $\pi$ satisfies (3).  In order to do so, we will assume now that $\pi$ is in nonincreasing order  $d_1 \geq \cdots \geq d_n.$

Since $d_n = k$ and $D = n-k-1,$ in $\pi,$ we need only show by (3) that

\begin{equation}
	\sum_{i=1}^{j} d_i \leq j(j-1) + \sum_{i=j+1}^{n} \text{min} \{j, d_i\}, \text{ for } k < j \leq n-k-1.
	\end{equation}

The sum on the right side of (4) is

$$\sum_{i=j+1}^{n} \text{ min } \{j,d_i\} = \sum_{i=j+1}^{n-k} \text{ min } \{j,d_i\} + \sum_{i=n-k+1}^{n} \text{ min } \{j,d_i\}.$$

Since $k <j \leq n-k-1,$ we have   $\text{ min } \{j,d_i\} = j$ for $i = j+1 \text{ to } n-k,$ while $\text{ min } \{j,d_i\} = k$ 
for $i =n-k+1 \text{ to } n.$  Therefore, 

$$\sum_{i=j+1}^{n} \text{ min } \{j,d_i\} = ((n-k) - (j+1) + 1)j + (n-(n-k+1)+1)k = (n-k-j)j  + k^2.$$

Thus, proving (4) reduces to showing 

\begin{equation}
	\sum_{i=1}^{j} d_i \leq j(j-1) + (n-k-j)j + k^2 = (n-k-1)j + k^2.
	\end{equation}

However,  $\sum_{i=1}^{j} d_i $ equals the sum $S_1$ of the largest j terms in the foundational $(n,k)-$sequence used to construct $\pi, $ plus  the sum $S_2 = \sum_{i=1}^{k-1} d_i' .$  Clearly, $S_1 \leq j(n-k-1) +1.$   Each term $d_i'$ in $S_2$ is at most $ k-1,$ and so $S_2 \leq (k-1)^2.$  Consequently, 

$$\sum_{i=1}^{j} d_i  = S_1 + S_2  \leq  j(n-k-1) +1  +(k-1)^2 < j(n-k-1) + k^2, $$

and so (5) holds and $\pi$ is forcibly hamiltonian.  Moreover, $\pi$ does not satisfy  Chv\'{a}tal's condition (1), and so $\pi$ is a Nash-Williams $(n,k)-$sequence.

Therefore,  given integers $n \geq 5 $  and $k \geq 2$, with $k < n/2,$ every sequence $\pi' = (d_1' \leq \cdots \leq d_{k-1}')$ of integers with  $ 0 \leq d_j' \leq k-1$ and $\pi' \neq (0, \cdots , 0, k-1)$ gives rise to a Nash-Williams $(n,k)-$sequence $\pi.$  Furthermore, we note that the map $\pi' \rightarrow  \pi$ is at most 2-to-1 and that the number of  $\pi'$ under consideration is $\binom{2(k-1)}{k-1}-1,$ and so the number of Nash-Williams $(n,k)-$sequences is at least $\frac{1}{2}[\binom{2(k-1)}{k-1}-1].$

We conclude that the number of Nash-Williams $n-$sequences, that is, $n-$seuqences that are forcibly hamiltonian by Theorem 2.3 but fail Chv\'{a}tal's condition (1) for some $k$, is at least

$$\sum_{k=2}^{\lfloor \frac{n-1}{2} \rfloor} \frac{1}{2}\Bigl[\binom{2(k-1)}{k-1}-1\Bigr] \geq \frac{1}{2} \Bigl[\binom{2(\lfloor \frac{n-1}{2} \rfloor-1)}{\lfloor \frac {n-1}{2} \rfloor -1} - 1\Bigr] =  \Omega ( \frac{2^n}{\sqrt n} ),$$

where the last equality holds by Stirling's formula.

\end{document}